\title{\bf{Some examples of algebraic surfaces with canonical map of degree 20}}
\author{
	NGUYEN BIN\\
	\\
	Dedicated to Margarida Mendes Lopes\\
	on the occasion of her sixty-fifth birthday
}
\date{\today}

\newcommand{\Addresses}{{
		\bigskip
		\footnotesize
			\text{Mathematics Division,}\par\nopagebreak	
			\text{National Center for Theoretical Sciences,}\par\nopagebreak	
			\text{Taiwan.}\par\nopagebreak		
		\textit{E-mail address}: \texttt{nguyenbin@ncts.ntu.edu.tw}
				
	}}

\newcommand\blfootnote[1]{%
	\begingroup
	\renewcommand\thefootnote{}\footnote{#1}%
	\addtocounter{footnote}{-1}%
	\endgroup
}

\date{\today}

\documentclass[10pt]{article}
\usepackage{amsmath, amsthm, amssymb, mathrsfs, amscd, amsthm, amscd,amsfonts,calligra,mathrsfs,adjustbox,lipsum}

\usepackage{indentfirst,url,xypic}
\usepackage[all]{xy}
\usepackage{url}
\usepackage{tikz}
\usepackage[colorlinks,plainpages]{hyperref}
\hypersetup{
	colorlinks=true,
	linkcolor=blue,
	filecolor=magenta,      
	urlcolor=cyan,
}

\DeclareMathOperator{\Picard}{Pic}

\newtheorem{Theorem}{Theorem}
\newtheorem{Proposition}{Proposition }

\newtheorem{Notation}[Proposition]{Notation}

\theoremstyle{remark}
\newtheorem{Remark}{Remark}

\newcommand{\MSC}{\textbf{Mathematics Subject Classification (2010):}}

\newcommand{\Key}{\textbf{Key words:}}

\begin{document}
\maketitle
\begin{abstract}  
	In this note, we construct two minimal surfaces of general type with geometric genus $ p_g = 3 $, irregularity $ q = 0 $, self-intersection of the canonical divisor $ K^2 = 20, 24 $ such that their canonical map is of degree $ 20 $. In one of these surfaces, the canonical linear system has a non-trivial fixed part. These surfaces, to our knowledge, are the first examples of minimal surfaces of general type with canonical map of degree $ 20 $.
\end{abstract}

\blfootnote{\MSC{ 14J29}.}
\blfootnote{\Key{ Surfaces of general type, Canonical maps, Abelian covers.}}

\section{Introduction} 
 If $ X $ is a minimal smooth complex projective surface, we denote by $ \xymatrix{\varphi_{\left| K_X\right| }:X \ar@{.>}[r] & \mathbb{P}^{p_g\left( X\right)-1}} $ the canonical map of $ X $, where $ K_{X} $ is the canonical divisor of $ X $ and $ p_g\left( X\right) = \dim H^0\left( X, K_{X}\right)  $ is the geometric genus. It is interesting to know which positive integers $ d $ occur as the degree of such
 canonical maps for surfaces of general type. This problem is motivated by the work of A. Beauville \cite{MR553705}. One knows that, for surfaces of general type, the degree $ d $ of the canonical map is at most $ 36 $ {\cite[ \rm Proposition 5.7]{MR527234}}. While surfaces with $ d = 1,2,3,\ldots,8 $ are easy to construct, only few surfaces with $ d > 8$ have been known so far. The first example was found by U. Persson \cite{MR527234} in 1977; in this example, the canonical map has degree $ 16 $. Then, a surface with $ d = 9 $ was constructed by S. L. Tan \cite{MR1141782} in 1992. In the last decade, some surfaces with $ d = 12,16,24,27,32,36 $ were constructed by C. Rito \cite{MR3391024}, \cite{MR3619737}, \cite{MR3663791}, \cite{2019arXiv190303017R}, C. Gleissner, R. Pignatelli and C. Rito \cite{2018arXiv180711854G}, Ching-Jui Lai and Sai-Kee Yeung \cite{2015arXiv151007097L}, and the author \cite{MR4008073}. In this paper, we present a way to construct surfaces with $ d = 20 $ as $ \mathbb{Z}_2^4 $-covers of the Del Pezzo surface $ Y_4 $ of degree $ 5 $.\\ 
 
 Throughout this paper all surfaces are projective algebraic over the complex numbers. The linear equivalence of divisors is denoted by $ \equiv $. We call a surface $ X $ no non-trivial $ 2 $-torsion if the only $ 2 $-torsion in $\Picard \left( X\right)  $ is $ \mathcal{O}_{X} $. A character $ \chi $ of the group $ \mathbb{Z}_{2}^4$ is a homomorphism from $ \mathbb{Z}_{2}^4$ to $ \mathbb{C}^{*} $, the multiplicative group of the non-zero complex numbers. We also use the following notations for Del Pezzo surfaces of degree $ 5 $:	
 \begin{Notation}\label{Notation of del Pezzo surface of degree 5}
 	We denote by $ Y_4 $ the blow-up of $ \mathbb{P}^2$ at four points in general position $ P_1, P_2, P_3, P_4 $. Let us denote by $ l $ the pull-back of a general line in $ \mathbb{P}^2$, by $ e_1 $, $ e_2 $, $ e_3 $, $ e_4 $ the exceptional divisors corresponding to $ P_1 $, $ P_2 $, $ P_3 $, $ P_4 $, respectively, by $ f_1 $, $ f_2$, $ f_3$, $ f_4$ the strict transforms of a general line through $ P_1 $, $ P_2 $, $ P_3 $, $ P_4 $, respectively and by $ h_{ij} $ the strict transforms of the line $ P_i  P_j $, for all $ i \neq j $ in $ \left\lbrace 1,2,3,4\right\rbrace  $, respectively. The anti-canonical class 
 	\begin{align*}
 		-K_{Y_4} &\equiv f_1 + f_2 + f_3 - e_4 \equiv f_1 + f_2 + f_4 - e_3 \equiv f_1 + f_3 + f_4 - e_2 \equiv f_2 + f_3 + f_4 - e_1  
 	\end{align*}
 	\noindent
 	is very ample and the linear system $ \left| -K_{Y_4} \right|  $ embeds $ Y_4 $ as a smooth Del Pezzo surface of degree $ 5 $ in $ \mathbb{P}^5 $.    	   		
 \end{Notation}
 
 The construction of abelian covers was studied by R. Pardini in \cite{MR1103912}. For details about the building data of abelian covers and their notations, we refer the reader to Section 1 and Section 2 of R. Pardini's work (\cite{MR1103912}). For the sake of completeness, we recall some facts on $ \mathbb{Z}_{2}^4 $-covers, in a form which is convenient for our later constructions. We will denote by  $ \chi_{j_1j_2j_3j_4} $ the character of $ \mathbb{Z}_{2}^4 $ defined by
 \begin{align*}
 	\chi_{j_1j_2j_3j_4}\left( a_1,a_2,a_3,a_4\right): =  e^{\left( \pi a_1j_1\right) \sqrt{-1}}e^{\left( \pi a_2j_2\right) \sqrt{-1}}e^{\left( \pi a_3j_3\right) \sqrt{-1}}e^{\left( \pi a_4j_4\right) \sqrt{-1}}
 \end{align*}
 for all $ j_1,j_2,j_3,j_4,a_1,a_2,a_3,a_4\in \mathbb{Z}_2 $. A $ \mathbb{Z}_{2}^4 $-cover $ \xymatrix{X \ar[r] & Y} $ can be determined by a collection of non-trivial divisors $ L_{\chi} $ labelled by characters of $ \mathbb{Z}_{2}^4 $ and effective divisors $ D_{\sigma} $ labelled by elements of $ \mathbb{Z}_{2}^4 $ of the surface $ Y $. More precisely, from \cite[\rm Theorem 2.1]{MR1103912} we can define $ \mathbb{Z}_{2}^4 $-covers as follows:
 \begin{Proposition} \label{Construction of cover of degree 8}
 	Given $ Y $ a smooth projective surface with no non-trivial $ 2 $-torsion, let $ L_{\chi} $ be divisors of $ Y $ such that $ L_{\chi} \not\equiv \mathcal{O}_Y $ for all non-trivial characters $ \chi $ of $ \mathbb{Z}_{2}^4  $ and let $ D_{\sigma} $ be effective divisors of  $ Y $ for all $ \sigma \in \mathbb{Z}_{2}^4 \setminus \left\lbrace \left(0,0,0,0 \right)  \right\rbrace  $ such that the total branch divisor $ B:=\sum\limits_{\sigma \ne 0}{D_{\sigma}} $ is reduced. Then $ \left\lbrace L_{\chi}, D_{\sigma} \right\rbrace_{\chi,\sigma}$ is the building data of a $ \mathbb{Z}_{2}^4$-cover $ \xymatrix{f:X \ar[r]& Y} $ if and only if
 	\begin{align}\label{The condition of Z_2^n covers} 
 		2L_{\chi} \equiv \sum\limits_{\chi\left( \sigma\right) = -1 }{D_{\sigma}}	
 	\end{align}
 	for all non-trivial characters $ \chi $ of $ \mathbb{Z}_{2}^4  $.
 	
 \end{Proposition}

 The following theorem is a result of this note: 
 \begin{Theorem}\label{the theorem with d = 20}
 	Let $ \xymatrix{f: X \ar[r]&Y_4} $ be a $ \mathbb{Z}_2^4 $-cover with the building data $ \left\lbrace L_{\chi}, D_{\sigma} \right\rbrace_{\chi,\sigma}$ such that the following hold:
 	\begin{enumerate}
 		\item \label{smooth condition}Each branch component $D_\sigma$ is smooth, the total branch locus $B $ is a simple normal crossings divisor and no more than two of these divisors $ D_{\sigma} $ go through the same point;
 		\item \label{factorization condition}$ D_{0100}+ D_{0101}+ D_{0110}+ D_{0111}$, $ D_{1000}+ D_{1001}+ D_{1010}+ D_{1011}$, $ D_{1100}+ D_{1101} + D_{1110}+ D_{1111}\in \left| -K_{Y_4}\right|$;
 		\item \label{generator condition}$ h^{0}\left( K_{Y_4} + L_{\chi} \right) =0$ for all $ \chi \notin \left\lbrace \chi_{1000}, \chi_{0100}, \chi_{1100}\right\rbrace  $;
 		\item \label{nef and big conditions}The divisor $ D_{0001}+ D_{0010}+ D_{0011} - K_{Y_4} $	is nef and big.	
 	\end{enumerate}
 	Then $ X $ is a minimal surface of general type with canonical map of degree $ 20 $ satisfying the following:
 	\begin{align*}
 		p_g\left( X\right) = 3, \hskip 0.3cm K_X^2 = 4\left( D_{0001}+ D_{0010}+ D_{0011} - K_{Y_4}\right)^2.
 	\end{align*}
 	Moreover, the reduced divisor supported on $  f^{*}\left( D_{0001}+ D_{0010}+ D_{0011}\right) $ is the fixed part of the canonical system $ \left| K_X \right| $.
 \end{Theorem} 
  
  \noindent
  Let us summarize the proof of Theorem \ref{the theorem with d = 20}. Assumptions (\ref{smooth condition}), (\ref{factorization condition}) and (\ref{nef and big conditions}) show that the surface $ X $ is a minimal surface of general type. Assumption (\ref{generator condition}) implies that the following diagram commutes (see Remark \ref{factorization of the main theorem} for the proof):
  $$
  
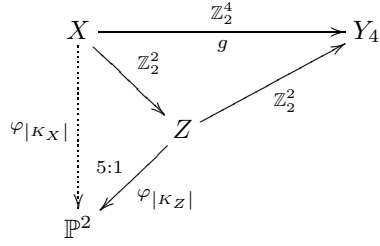
\xymatrix{X \ar[0,3]^{\mathbb{Z}_2^4}_g \ar[1,1]^{\mathbb{Z}_2^2} \ar@{.>}[2,0]_{\varphi_{\left| K_X \right| }}&&& Y_4\\
  	&Z \ar[-1,2]_{\mathbb{Z}_2^2} \ar[1,-1]^{\varphi_{\left| K_{Z} \right| }}_{5:1}&&\\
  	\mathbb{P}^{2}&&&} 	
  $$
  
   \noindent
  In the above diagram, the intermediate surface $ Z:=X/\Gamma$ is the quotient surface of $ X$, where $ \Gamma:= \left\langle \left( 0,0,0,1\right), \left( 0,0,1,0\right) \right\rangle  $ is the subgroup of $ \mathbb{Z}_2^4 $. The surface $ Z $ is the bidouble cover of $ Y_4 $ ramified on 
  \begin{align*}
  \left( D_{0100}+ D_{0101}+ D_{0110}+ D_{0111}\right) +\left(  D_{1000}+ D_{1001}+ D_{1010}+ D_{1011}\right) +\left(  D_{1100}+ D_{1101} + D_{1110}+ D_{1111}\right).	
  \end{align*}
  
  \noindent
  Assumption (\ref{factorization condition}) shows that the canonical map of $ Z $ is of degree $ 5 $ (see Remark \ref{factorization of the main theorem} for the proof). Therefore, the canonical map of $ X $ is of degree $ 20 $. As application of Theorem \ref{the theorem with d = 20}, we construct two surfaces with $ d =20 $ described as follows: 
  
  \begin{Theorem}\label{the main theorem with d = 20}
  	There exist minimal surfaces of general type $ X $ satisfying the following
  	$$
  	\begin{tabular}{|c |c| c| c| c|}
  		\hline
  		$ d $&$ K_X^2 $ &$ p_g\left( X\right) $ &$ q\left( X\right) $& $\left| K_X \right| $ \\
  		\hline
  		$ 20 $&$ 20 $&$ 3 $&$ 0 $& base point free\\
  		\hline
  		$ 20 $&$ 24 $&$ 3 $&$ 0 $& has a non-trivial fixed part\\
  		\hline
  	\end{tabular} 
  	$$
  \end{Theorem}

	\section{$ \mathbb{Z}_{2}^4 $-coverings} 
	
	For the convenience of the reader, we leave here the relations (\ref{The condition of Z_2^n covers}) of the building data of $ \mathbb{Z}_{2}^4 $-covers:
	$$
	\begin{adjustbox}{max width=\textwidth}
	\begin{tabular}{l l r r r r r r r r r r r r r r}
	\color{red}$ B $&\color{red}$ =D_{0001 } $&\color{red}$ +D_{0010 } $&\color{red}$ +D_{0011} $&\color{red}$ +D_{0100} $&\color{red}$ +D_{0101} $&\color{red}$ +D_{0110} $&\color{red}$ +D_{0111} $ &\color{red}$ +D_{1000} $ &\color{red}$ +D_{1001} $ &\color{red}$ +D_{1010} $ &\color{red}$ +D_{1011} $ &\color{red}$ +D_{1100} $ &\color{red}$ +D_{1101} $ &\color{red}$ +D_{1110} $ &\color{red}$ +D_{1111} $  \\
	\color{blue}$ 2L_{0001} $&\color{blue}$ \equiv D_{0001 } $&\color{blue}$  $&\color{blue}$ +D_{0011} $&\color{blue}$  $&\color{blue}$ +D_{0101} $&\color{blue}$  $&\color{blue}$ +D_{0111} $ &\color{blue}$  $ &\color{blue}$ +D_{1001} $ &\color{blue}$  $ &\color{blue}$ +D_{1011} $ &\color{blue}$  $ &\color{blue}$ +D_{1101} $ &\color{blue}$  $ &\color{blue}$ +D_{1111} $ \\ 
	\color{blue}$ 2L_{0010} $&\color{blue}$ \equiv $&\color{blue}$ D_{0010 } $&\color{blue}$ +D_{0011} $&\color{blue}$ $&\color{blue}$  $&\color{blue}$ +D_{0110} $&\color{blue}$ +D_{0111} $ &\color{blue}$  $ &\color{blue}$  $ &\color{blue}$ +D_{1010} $ &\color{blue}$ +D_{1011} $ &\color{blue}$  $ &\color{blue}$  $ &\color{blue}$ +D_{1110} $ &\color{blue}$ +D_{1111} $ \\ 
	\color{blue}$ 2L_{0100} $&\color{blue}$ \equiv $&\color{blue}$  $&\color{blue}$  $&\color{blue}$ D_{0100} $&\color{blue}$ +D_{0101} $&\color{blue}$ +D_{0110} $&\color{blue}$ +D_{0111} $ &\color{blue}$  $ &\color{blue}$  $ &\color{blue}$  $ &\color{blue}$  $ &\color{blue}$ +D_{1100} $ &\color{blue}$ +D_{1101} $ &\color{blue}$ +D_{1110} $ &\color{blue}$ +D_{1111} $ \\
	\color{blue}$ 2L_{1000} $&\color{blue}$ \equiv $&\color{blue}$  $&\color{blue}$  $&\color{blue}$  $&\color{blue}$  $&\color{blue}$  $&\color{blue}$  $ &\color{blue}$ D_{1000} $ &\color{blue}$ +D_{1001} $ &\color{blue}$ +D_{1010} $ &\color{blue}$ +D_{1011} $ &\color{blue}$ +D_{1100} $ &\color{blue}$ +D_{1101} $ &\color{blue}$ +D_{1110} $ &\color{blue}$ +D_{1111} $ \\
	\color{black}$ 2L_{0011} $&$ \equiv D_{0001 } $&$ +D_{0010 } $&$ $&$  $&$ +D_{0101} $&$ +D_{0110} $&$ $ &$  $ &$ +D_{1001} $ &$ +D_{1010} $ &$  $ &$  $ &$ +D_{1101} $ &$ +D_{1110} $ &$  $ \\
	$ 2L_{0101} $&$ \equiv D_{0001 } $&$ $&$ +D_{0011} $&$ +D_{0100} $&$  $&$ +D_{0110} $&$  $ &$  $ &$ +D_{1001} $ &$  $ &$ +D_{1011} $ &$ +D_{1100} $ &$  $ &$ +D_{1110} $ &$  $ \\
	$ 2L_{0110} $&$ \equiv $&$ D_{0010 } $&$ +D_{0011} $&$ +D_{0100} $&$ +D_{0101} $&$  $&$  $ &$  $ &$  $ &$ +D_{1010} $ &$ +D_{1011} $ &$ +D_{1100} $ &$ +D_{1101} $ &$  $ &$  $ \\
	$ 2L_{0111} $&$ \equiv D_{0001 } $&$ +D_{0010 } $&$ $&$ +D_{0100} $&$  $&$  $&$ +D_{0111} $ &$  $ &$ +D_{1001} $ &$ +D_{1010} $ &$  $ &$ +D_{1100} $ &$  $ &$  $ &$ +D_{1111} $ \\
	$ 2L_{1001} $&$ \equiv D_{0001 } $&$  $&$ +D_{0011} $&$  $&$ +D_{0101} $&$  $&$ +D_{0111} $ &$ +D_{1000} $ &$  $ &$ +D_{1010} $ &$  $ &$ +D_{1100} $ &$  $ &$ +D_{1110} $ &$  $ \\
	$ 2L_{1010} $&$ \equiv $&$ D_{0010 } $&$ +D_{0011} $&$  $&$  $&$ +D_{0110} $&$ +D_{0111} $ &$ +D_{1000} $ &$ +D_{1001} $ &$  $ &$  $ &$ +D_{1100} $ &$ +D_{1101} $ &$  $ &$ $ \\
	$ 2L_{1011} $&$ \equiv D_{0001 } $&$ +D_{0010 } $&$ $&$  $&$ +D_{0101} $&$ +D_{0110} $&$ $ &$ +D_{1000} $ &$$ &$  $ &$ +D_{1011} $ &$ +D_{1100} $ &$  $ &$  $ &$ +D_{1111} $ \\
	$ 2L_{1100} $&$ \equiv $&$  $&$  $&$ D_{0100} $&$ +D_{0101} $&$ +D_{0110} $&$ +D_{0111} $ &$ +D_{1000} $ &$ +D_{1001} $ &$ +D_{1010} $ &$ +D_{1011} $ &$  $ &$  $ &$  $ &$  $ \\
	$ 2L_{1101} $&$ \equiv D_{0001 } $&$  $&$ +D_{0011} $&$ +D_{0100} $&$  $&$ +D_{0110} $&$  $ &$ +D_{1000} $ &$  $ &$ +D_{1010} $ &$  $ &$  $ &$ +D_{1101} $ &$  $ &$ +D_{1111} $ \\
	$ 2L_{1110} $&$ \equiv $&$ D_{0010 } $&$ +D_{0011} $&$ +D_{0100} $&$ +D_{0101} $&$  $&$  $ &$ +D_{1000} $ &$ +D_{1001} $ &$  $ &$  $ &$  $ &$  $ &$ +D_{1110} $ &$ +D_{1111} $ \\
	$ 2L_{1111} $&$ \equiv D_{0001 } $&$ +D_{0010 } $&$  $&$ +D_{0100} $&$  $&$  $&$ +D_{0111} $ &$ +D_{1000} $ &$  $ &$  $ &$ +D_{1011} $ &$  $ &$ +D_{1101} $ &$ +D_{1110} $ &$  $ \\
\end{tabular}
\end{adjustbox}
$$
	
	By \cite[\rm Theorem 3.1]{MR1103912} if each branch component $D_\sigma$ is smooth and the total branch locus $B $ is a simple normal crossings divisor, the surface $X$ is smooth. \\
	
	Also from \cite[\rm Lemma 4.2, Proposition 4.2]{MR1103912} we have:
	\begin{Proposition}\label{invariants of Z_2^n cover}
	If $ Y $ is a smooth surface and $ \xymatrix{f: X \ar[r]& Y} $ is a smooth $  \mathbb{Z}_{2}^4$-cover with the building data $ \left\lbrace L_{\chi}, D_{\sigma} \right\rbrace_{\chi,\sigma}$, the surface $ X $ satisfies the following:
		\begin{align}
		2K_X & \equiv f^*\left( 2K_Y + \sum\limits_{\sigma \ne 0} {D_{\sigma} } \right); \\
		f_{*}\mathcal{O}_X &= \mathcal{O}_Y \oplus \bigoplus\limits_{\chi \ne \chi_{0000}  }L_{\chi}^{-1};\\	
	    H^{0}\left( X, K_X\right) &= H^{0}\left( Y, K_{Y}\right) \oplus \bigoplus_{\chi \ne  \chi_{0000}}{H^{0}\left( Y, K_{Y} +L_{\chi}\right)}; \label{property1}\\
	    K^2_X &= 4\left( 2K_Y + \sum\limits_{\sigma \ne 0} {D_{\sigma} } \right)^2; \\
	    p_g\left( X \right) &=p_g\left( Y \right) +\sum\limits_{\chi \ne  \chi_{0000}  }{h^0\left( L_{\chi} + K_Y \right)}; \\
	    \chi\left( \mathcal{O}_X \right) &= 16\chi\left( \mathcal{O}_Y \right)  +\sum\limits_{\chi \ne \chi_{0000}  }{\frac{1}{2}L_{\chi}\left( L_{\chi}+K_Y\right)}. 
	\end{align}
    \noindent
    Moreover, the canonical linear system $ \left|  K_X \right|  $ is generated by
    \begin{align}
    	f^{*}\left(  \left| K_Y + L_{\chi}\right| \right) +\sum\limits_{\chi\left( \sigma\right)=1 }{R_{\sigma}}    , \hskip 5pt \forall \chi \in J  \label{The generators of canonical system of Z_2^n cover}
    \end{align}
    \noindent
    where $ J:= \left\lbrace  \chi' : \left| K_Y + L_{\chi'}\right| \ne \emptyset \right\rbrace  $ and $ R_{\sigma} $ is the reduced divisor supported on $ f^{*}\left( D_{\sigma} \right)  $.
	\end{Proposition}

    \noindent
    For the proof of the last statement of Proposition \ref{invariants of Z_2^n cover}, we refer the reader to \cite[\rm Page 3]{2018arXiv180711854G}.	
   
    \section{Surfaces with $ d = 20 $ as $ \mathbb{Z}_2^4 $-covers}
    \subsection{Proof of Theorem \ref{the theorem with d = 20}}
    The surface $ X $ is smooth because each branch component $ D_{\sigma} $ is smooth, the total branch locus $ B $ is a normal crossings divisor and no more than two of these divisors $ D_{\sigma} $ go through the same point. Moreover, by Proposition \ref{invariants of Z_2^n cover}, the surface $ X $ satisfies the following:		
    \begin{align*}
    	2K_X &\equiv f^*\left( 2K_{Y_4} + \sum\limits_{\sigma}{D_{\sigma}}\right)\\
    	&\equiv f^*\left( D_{0001}+ D_{0010}+ D_{0011} -K_{Y_4}  \right).
    \end{align*}
    \noindent
    We notice that a surface is of general type and minimal if the canonical divisor is big and nef (see e.g. \cite[\rm Section 2]{MR2931875}). We remark that the divisor $ D_{0001}+ D_{0010}+ D_{0011} -K_{Y_4} $ is nef and big by Assumption (\ref{nef and big conditions}). Since the divisor $ 2K_X $ is the pull-back of a nef and big divisor, the canonical divisor $ K_X $ is nef and big. Thus, the surface $ X $ is of general type and minimal. Furthermore, from Proposition \ref{invariants of Z_2^n cover}, the surface $X$ possesses the following invariants:
    \begin{align*}
    	p_g\left( X\right) = 3, \hskip 0.3cm K_X^2= 4\left( D_{0001}+ D_{0010}+ D_{0011} -K_{Y_4}  \right)^2. 
    \end{align*}
    
    We show that the canonical map $ \varphi_{\left| K_X  \right|}  $ has degree $ 20 $. By Assumptions (\ref{factorization condition}) and (\ref{generator condition}), we have 
    \begin{align*}
    	&L_{1000} + K_{Y_4} \equiv 	L_{0100} + K_{Y_4} \equiv 	L_{1100} + K_{Y_4} \equiv \mathcal{O}_{Y_4},\\
    	&h^{0}\left(  L_{\chi} +K_{Y_4}\right) =0, \hskip0.5cm  \forall \chi \notin \left\lbrace \chi_{1000}, \chi_{0100}, \chi_{1100}\right\rbrace.
    \end{align*}
    \noindent
    By (\ref{The generators of canonical system of Z_2^n cover}), the linear system $ \left| K_X\right|  $ is generated by the three following divisors:
    \begin{align*}
    	\overline{D}_{0001}+\overline{D}_{0010}+\overline{D}_{0011} + \overline{D}_{0100} + \overline{D}_{0101} + \overline{D}_{0110} + \overline{D}_{0111},\\
    	\overline{D}_{0001}+\overline{D}_{0010}+\overline{D}_{0011} + \overline{D}_{1000} + \overline{D}_{1001} + \overline{D}_{1010} + \overline{D}_{1011},\\
    	\overline{D}_{0001}+\overline{D}_{0010}+\overline{D}_{0011} + \overline{D}_{1100} + \overline{D}_{1101} + \overline{D}_{1110} + \overline{D}_{1111},
    \end{align*}
    \noindent
    where $ \overline{D}_{\sigma}$ are the reduced divisors supported $f^{*}\left( D_{\sigma}\right)  $, for all $ \sigma $. Because the divisors $ \overline{D}_{0001}$, $\overline{D}_{0010}$, $\overline{D}_{0011} $ are common components of the three above divisors, these divisors $ \overline{D}_{0001}$, $\overline{D}_{0010}$, $\overline{D}_{0011} $ are fixed components of $ \left| K_X\right|  $.
    
    \noindent
    On the other hand, by Assumption (\ref{smooth condition}) the three divisors $ \overline{D}_{0100} + \overline{D}_{0101} + \overline{D}_{0110} + \overline{D}_{0111} $, $ \overline{D}_{1000} + \overline{D}_{1001} + \overline{D}_{1010} + \overline{D}_{1011} $, $ \overline{D}_{1100} + \overline{D}_{1101} + \overline{D}_{1110} + \overline{D}_{1111} $ have no common intersection. So the linear system $ \left| M\right|  $ is base point free, where $ M:= \overline{D}_{0100} + \overline{D}_{0101} + \overline{D}_{0110} + \overline{D}_{0111}$. This together with $ M^2 =4\left( 3l -e_1 - e_2 - e_3 -e_4\right)^2 = 20 > 0 $ implies that the linear system $ \left| K_X\right|  $ is not composed with a pencil. Thus, the canonical image is $ \mathbb{P}^2 $, the canonical map is of degree $ 20 $, and the divisor $ \overline{D}_{0001}+\overline{D}_{0010}+\overline{D}_{0011} $ is the fixed part of $ \left| K_X\right|  $.
    
    \begin{Remark}\label{factorization of the main theorem}
    The canonical map $ \varphi_{\left| K_X  \right|} $ of $ X $ is the composition of the quotient map $ \xymatrix{X \ar[r]& Z:= X/\Gamma } $ with the canonical map $\varphi_{\left| K_{Z} \right| }$ of $ Z $. Moreover, the canonical map of $ Z $ is of degree $ 5 $.	
    \end{Remark}
    \noindent
    In fact, by  (\ref{property1}), we have the following decomposition:     
    \begin{align*}
    	H^{0}\left( X, K_X\right) = H^{0}\left( Y_4, K_{Y_4}\right) \oplus \bigoplus_{\chi \ne  \chi_{0000}}{H^{0}\left( Y_4, K_{Y_4} +L_{\chi}\right)}. 
    \end{align*}
    
    \noindent
    The group $ \Gamma:= \left\langle \left( 0,0,0,1\right), \left( 0,0,1,0\right) \right\rangle  $ is the subgroup of $ \mathbb{Z}_2^4 $. Let $ \Gamma^\perp $ denote the kernel of the restriction map $ \xymatrix{\left( \mathbb{Z}_2^4\right)^{*} \ar[r]&\Gamma^{*}} $, where $ \Gamma^{*} $ is the character group of $ \Gamma $. We have $ \Gamma^{\perp} = \left\langle \chi_{1000},\chi_{0100},\chi_{1100} \right\rangle  $. The subgroup $ \Gamma $ acts trivially on $ H^{0}\left( X, K_X\right) $ since $ h^0\left( L_{ \chi} + K_{Y_4} \right)  = 0 $ for all $ \chi \notin \Gamma^\perp $ by Assumption (\ref{generator condition}). So the canonical map $ \varphi_{\left| K_X \right| } $ is the composition of the quotient map $ \xymatrix{X \ar[r]& Z:= X/\Gamma } $ with the canonical map $\varphi_{\left| K_{Z} \right| }$ of $ Z $ (see e.g. \cite[\rm Example 2.1]{MR1103913}). 
    
    \noindent
    The intermediate surface $ Z $ is the bidouble cover of $ Y_4 $ with the building data $ \left\lbrace D_1, D_2, D_3, L_1, L_2, L_3\right\rbrace  $ determined as follows:
    \begin{align*}
    	D_1:= &D_{0100} + D_{0101} +D_{0110} + D_{0111} \equiv -K_{Y_4}, &L_1:= &L_{1000} \equiv -K_{Y_4},\\ 
    	D_2:= &D_{1000} + D_{1001} +D_{1010} + D_{1011} \equiv -K_{Y_4}, &L_2:= &L_{0100} \equiv -K_{Y_4},\\
    	D_3:= &D_{1100} + D_{1101} +D_{1110} + D_{1111} \equiv -K_{Y_4}, &L_3:= &L_{1100} \equiv -K_{Y_4}.
    \end{align*}
    \noindent
    Assumption (\ref{smooth condition}) shows that the singularities of $ Z $ are nodes and the canonical map of $ Z $ is of degree $ \left( 3l - e_1 - e_2 - e_3 - e_4\right)^2 = 5 $.  	
   
    \subsection{Constructions of the surfaces in Theorem \ref{the main theorem with d = 20}}

	\subsubsection{A surface with $ d = 20 $, $ p_g = 3 $, $ q = 0 $, $ K^2 = 20 $}
	In this section, we construct the surface described in the first row of Theorem \ref{the main theorem with d = 20}. Let $ Y_4 $ be a Del Pezzo surface of degree $ 5 $ (see Notation \ref{Notation of del Pezzo surface of degree 5}). We consider the following smooth divisors of $ Y_4 $:
	\begin{align*}
	D_{0101}:=& h_{14}         & D_{0110}:=&f_{31}+e_{1}& D_{0111}:=&h_{12} \\
	D_{1001}:=& f_{11}+e_{2}   & D_{1010}:=&h_{23}      & D_{1011}:=&h_{24} \\
	D_{1101}:=& h_{13}         & D_{1110}:=&h_{34}      & D_{1111}:=&f_{21}+e_{3}  
	\end{align*}
	\noindent
	and $ D_{\sigma} = 0 $ for the other $ \sigma $, where $ f_{11} \in \left| f_1\right| $, $ f_{21} \in \left| f_2\right| $ and $ f_{31} \in \left| f_3\right| $ such that no more than two of these divisors $ D_{\sigma} $ go through the same point. We consider the following non-trivial divisors of $ Y_4 $: 
	$$
	\begin{tabular}{l r r r r}
	$ L_{0001}:= $&$ 2f_{1} $&$+f_{2} $&$ $&$ -e_{4}$ \\
	$ L_{0010}:= $&$  $&$ 2f_{2} $&$+f_{3} $&$ -e_{4}$ \\
	$ L_{0100}:= $&$ f_{1} $&$ +f_{2} $&$+f_{3} $&$ -e_{4}$ \\
	$ L_{1000}:= $&$ f_{1} $&$ +f_{2} $&$+f_{3} $&$ -e_{4}$\\
	$ L_{0011}:= $&$ f_{1}$&$ $&$ +2f_{3}$&$ -e_{4}$\\
	$ L_{0101}:= $&$   $&$ $&$ f_{3}$&$ +f_{4}$\\
	$ L_{0110}:= $&$ h_{12}  $&$ +h_{34} $&$ $&$ $\\
	$ L_{0111}:= $&$  f_{1}  $&$+f_{2} $&$ $&$ $ \\
	$ L_{1001}:= $&$ h_{12}  $&$ +h_{34} $&$ $&$ $\\
	$ L_{1010}:= $&$  f_{1}  $&$ $&$ +f_{3}$&$ $ \\
	$ L_{1011}:= $&$   $&$ f_{2} $&$ $&$ +f_{4}$\\
	$ L_{1100}:= $&$ f_{1} $&$ +f_{2} $&$+f_{3} $&$-e_{4} $\\
	$ L_{1101}:= $&$   $&$ f_{2} $&$ +f_{3}$&$ $\\
	$ L_{1110}:= $&$ f_{1} $&$  $&$ $&$+f_{4} $ \\
	$ L_{1111}:= $&$ h_{12}  $&$ +h_{34} $.&$ $&$ $
	\end{tabular}
	$$
	\noindent
	These divisors $ D_{\sigma}, L_{\chi} $ satisfy the following relations:
	$$
	\begin{adjustbox}{max width=\textwidth}
	\begin{tabular}{l l l l l l l l l l l r r r r r}
	$ 2L_{0001} $& $ \equiv $& $ D_{0101} $ & $  $ &$ +D_{0111} $& $ +D_{1001} $& $  $ &$ +D_{1011} $& $ +D_{1101} $& $  $ &$ +D_{1111} $&$ \equiv $&$ 4f_{1} $&$+2f_{2} $&$ $&$ -2e_{4}$\\
	$ 2L_{0010} $& $ \equiv $& $  $ & $  D_{0110} $ &$ +D_{0111} $& $  $& $ +D_{1010} $ &$ +D_{1011} $& $  $& $ +D_{1110} $& $ +D_{1111} $&$ \equiv $&$  $&$ 4f_{2} $&$+2f_{3} $&$ -2e_{4}$\\
	$ 2L_{0100} $& $ \equiv $& $ D_{0101} $ & $ + D_{0110} $ &$ +D_{0111} $& $  $& $  $ &$  $& $ +D_{1101} $& $ +D_{1110} $& $ +D_{1111} $&$ \equiv $&$ 2f_{1} $&$ +2f_{2} $&$+2f_{3} $&$ -2e_{4}$\\
	$ 2L_{1000} $& $ \equiv $& $  $ & $  $ &$  $& $ D_{1001} $& $ +D_{1010} $ &$ +D_{1011} $& $ +D_{1101} $& $ +D_{1110} $& $ +D_{1111} $&$ \equiv $&$ 2f_{1} $&$ +2f_{2} $&$+2f_{3} $&$ -2e_{4}$\\
	$ 2L_{0011} $& $ \equiv $& $ D_{0101} $ & $ + D_{0110} $ &$  $& $ +D_{1001} $& $ +D_{1010} $ &$  $& $ +D_{1101} $& $ +D_{1110} $& $  $&$ \equiv $&$ 2f_{1}$&$ $&$ +4f_{3}$&$ -2e_{4}$\\
	$ 2L_{0101} $& $ \equiv $& $  $ & $  D_{0110} $ &$  $& $ +D_{1001} $& $  $ &$ +D_{1011} $& $  $& $ +D_{1110} $& $  $&$ \equiv $&$   $&$ $&$ 2f_{3}$&$ +2f_{4}$\\
	$ 2L_{0110} $& $ \equiv $& $ D_{0101} $ & $  $ &$  $& $  $& $ +D_{1010} $ &$ +D_{1011} $& $ +D_{1101} $& $  $& $  $&$ \equiv $&$ 2h_{12}  $&$ +2h_{34} $&$ $&$ $\\
	$ 2L_{0111} $& $ \equiv $& $  $ & $  $ &$ D_{0111} $& $ +D_{1001} $& $ +D_{1010} $ &$  $& $  $& $  $& $ +D_{1111} $&$ \equiv $&$  2f_{1}  $&$+2f_{2} $&$ $&$ $\\
	$ 2L_{1001} $& $ \equiv $& $ D_{0101} $ & $  $ &$ +D_{0111} $& $  $& $ +D_{1010} $ &$  $& $  $& $ +D_{1110} $& $  $&$ \equiv $&$ 2h_{12}  $&$ +2h_{34} $&$ $&$ $\\
	$ 2L_{1010} $& $ \equiv $& $  $ & $  D_{0110} $ &$ +D_{0111} $& $ +D_{1001} $& $  $ &$  $& $ +D_{1101} $& $  $& $  $&$ \equiv $&$  2f_{1}  $&$ $&$ +2f_{3}$&$ $\\
	$ 2L_{1011} $& $ \equiv $& $ D_{0101} $ & $ + D_{0110} $ &$  $& $  $& $  $ &$ +D_{1011} $& $  $& $  $& $ +D_{1111} $&$ \equiv $&$   $&$ 2f_{2} $&$ $&$ +2f_{4}$\\
	$ 2L_{1100} $& $ \equiv $& $ D_{0101} $ & $ + D_{0110} $ &$ +D_{0111} $& $ +D_{1001} $& $ +D_{1010} $ &$ +D_{1011} $& $  $& $  $& $  $&$ \equiv $&$ 2f_{1} $&$ +2f_{2} $&$+2f_{3} $&$-2e_{4} $\\
	$ 2L_{1101} $& $ \equiv $& $  $ & $  D_{0110} $ &$  $& $  $& $ +D_{1010} $ &$  $& $ +D_{1101} $& $  $& $ +D_{1111} $&$ \equiv $&$   $&$ 2f_{2} $&$ +2f_{3}$&$ $\\
	$ 2L_{1110} $& $ \equiv $& $ D_{0101} $ & $  $ &$  $& $ +D_{1001} $& $  $ &$  $& $  $& $ +D_{1110} $& $ +D_{1111} $&$ \equiv $&$ 2f_{1} $&$  $&$ $&$+2f_{4} $\\
	$ 2L_{1111} $& $ \equiv $& $  $ & $  $ &$ D_{0111} $& $  $& $  $ &$ +D_{1011} $& $ +D_{1101} $& $ +D_{1110} $& $  $&$ \equiv $&$ 2h_{12}  $&$ +2h_{34} $.&$ $&$ $
\end{tabular}
\end{adjustbox}
$$
	
	\noindent
	Thus by Proposition \ref{Construction of cover of degree 8}, the divisors $ D_{\sigma}, L_{\chi} $ define a $ \mathbb{Z}^4_2 $-cover $ \xymatrix{g: X \ar[r] & Y_4}  $. Moreover, this $ \mathbb{Z}^4_2 $-cover fulfils the hypotheses of Theorem \ref{the theorem with d = 20}. In fact, we have that
	\begin{align*}
		D_{0100} + D_{0101} +D_{0110} + D_{0111} &=h_{14}  + f_{31}+e_{1} + h_{12}\equiv 3l - e_1 - e_2 - e_3 - e_4\\ 
		D_{1000} + D_{1001} +D_{1010} + D_{1011} &=f_{11}+e_{2} + h_{23} + h_{24}\equiv 3l - e_1 - e_2 - e_3 - e_4\\
		D_{1100} + D_{1101} +D_{1110} + D_{1111} &=h_{13} + h_{34}  + f_{21}+e_{3}\equiv 3l - e_1 - e_2 - e_3 - e_4,
	\end{align*}
	\noindent
	$ h^{0}\left( K_{Y_4} + L_{\chi} \right) =0$ for all $ \chi \notin \left\lbrace \chi_{1000}, \chi_{0100}, \chi_{1100}\right\rbrace  $, and the divisor $ D_{0001}+ D_{0010}+ D_{0011} - K_{Y_4} \equiv 3l - e_1 - e_2 - e_3 - e_4$ is nef and big. Thus by Theorem \ref{the theorem with d = 20} and Proposition \ref{invariants of Z_2^n cover}, the surface $ X $ is a minimal surface of general type and possesses the following invariants:
	\begin{align*}
		K_X^2= 20, p_g\left( X\right) = 3, \chi\left( \mathcal{O}_X\right) =4, q\left( X\right)  = 0.   
	\end{align*}
	\noindent
	Moreover, the canonical map $ \varphi_{\left| K_X  \right|}  $ is of degree $ 20 $ and the linear system $ \left| K_X\right|  $ is base point free.\\
	
	\begin{Remark}
	The surface $ X $ has four pencils of genus $ 9 $ corresponding to the fibres $ f_1, f_2, f_3, f_4$.
	\end{Remark}

    \vskip 0.5cm
	In the above construction, for each choice of $ f_{11} \in \left| f_1\right| $, $ f_{21} \in \left| f_2\right| $ and $ f_{31} \in \left| f_3\right| $, we obtain a natural deformation of the surface $ X $ (we refer \cite[\rm Definition 5.1]{MR1103912} for the definition of natural deformations of an abelian cover). It is worth pointing out that a natural deformation of an abelian cover $ \xymatrix{X \ar[r]& Y} $ is a deformation of the map $ \xymatrix{X \ar[r]& Y} $ by \cite[\rm Proposition 5.1]{MR1103912}.
	\begin{Remark}\label{deformation of surface}
		The surface $ X $ admits natural deformations. Moreover, all the natural deformations of $ X $ are Galois.
	\end{Remark}
    \noindent
    In fact, by \cite[\rm Definition 5.1]{MR1103912} the natural deformations of the $ \mathbb{Z}_2^4 $-cover $ \xymatrix{g: X \ar[r] & Y_4}  $ are parametrized by the direct sum of the vector spaces
    \begin{align*}
    	\bigoplus \limits_{\sigma \neq 0  }{H^{0}\left( Y_4, D_{\sigma}\right) } \bigoplus \bigoplus \limits_{\substack{\sigma \neq 0 \\ \chi \ne \chi_{0000} \\\chi\left( \sigma\right) =1}  }{H^{0}\left( Y_4, D_{\sigma} - L_{\chi}\right) }.
    \end{align*}
    \noindent
    Moreover, all the natural deformations of $ X $ are Galois if the second summand $ \bigoplus \limits_{\substack{\sigma \neq 0 \\ \chi \ne \chi_{0000} \\\chi\left( \sigma\right) =1}  }{H^{0}\left( Y_4, D_{\sigma} - L_{\chi}\right) } $ is zero (see \cite[\rm Definition 3.2]{MR1444010}). We have that
    \begin{align*}
    	H^{0}\left( Y_4, D_{0110}\right) &= H^{0}\left( Y_4, f_{31}\right) \cong \mathbb{C}^2\\
    	H^{0}\left( Y_4, D_{1001}\right) &= H^{0}\left( Y_4, f_{11}\right) \cong \mathbb{C}^2\\
    	H^{0}\left( Y_4, D_{1111}\right) &= H^{0}\left( Y_4, f_{21}\right) \cong \mathbb{C}^2
    \end{align*}
    \noindent
    and $ H^{0}\left( Y_4, D_{\sigma}\right) \cong \mathbb{C} $ for the other non-trivial $ D_{\sigma} $. So the famify of natural deformations of $\xymatrix{g: X \ar[r] & Y_4} $ is parametrized by the base space $ \mathbb{P}^1\times \mathbb{P}^1\times\mathbb{P}^1$. Furthermore, all natural deformations of $ X $ are Galois since $ \bigoplus \limits_{\substack{\sigma \neq 0 \\ \chi \ne \chi_{0000} \\\chi\left( \sigma\right) =1}  }{H^{0}\left( Y_4, D_{\sigma} - L_{\chi}\right) } =0$.
	 	  
	  \subsubsection{A surface with $ d = 20 $, $ p_g = 3 $, $ q = 0 $, $ K^2 = 24 $}\label{The main construction 2 of surfaces with d = 20}
	  In this section, we construct the surface described in the second row of Theorem \ref{the main theorem with d = 20}. We consider the following smooth divisors of a del Pezzo surface $ Y_4 $ of degree $ 5 $:
	  \begin{align*}
	  D_{0011}:=& e_{4}         & &&  \\
	  D_{0101}:=& h_{14}  & D_{0110}:=&f_{21}          & D_{0111}:=&f_{31} & &\\
	  D_{1000}:=& e_{2}   & D_{1001}:=&h_{23}          & D_{1010}:=&h_{24} & D_{1011}:=&f_{11}\\
	  D_{1100}:=& h_{34}  & D_{1101}:=&h_{12} +h_{13}  & D_{1110}:=&e_{1}  & D_{1111}:=&e_{3} 
	  \end{align*}
	  \noindent
	  and the other $ D_{\sigma} = 0 $, where $ f_{11} \in \left| f_1\right| $, $ f_{21} \in \left| f_2\right| $ and $ f_{31} \in \left| f_3\right| $ such that no more than two of these divisors $ D_{\sigma} $ go through the same point. We consider the following non-trivial divisors of $ Y_4 $:
	  $$
	  \begin{tabular}{l r r r r r}
	  $ L_{0001}:= $&$ 2f_{1} $&$+f_{2} $&$ $&$-e_{3} $& \\
	  $ L_{0010}:= $&$  $&$ f_{2} $&$ $&$+l $& \\
	  $ L_{0100}:= $&$ f_{1} $&$ +f_{2} $&$+f_{3} $&$ $&$ -e_{4}$ \\
	  $ L_{1000}:= $&$ f_{1} $&$ +f_{2} $&$+f_{3} $&$ $&$ -e_{4}$\\
	  $ L_{0011}:= $&$ f_{1}$&$+2f_{2} $&$ $&$-e_{3} $&$ -e_{4}$\\
	  $ L_{0101}:= $&$   $&$f_{2} $&$ +f_{3}$&$ $&\\
	  $ L_{0110}:= $&$ 2f_{1}  $&$ +f_{2} $&$ $&$-e_{3}  $&$ -e_{4}$\\
	  $ L_{0111}:= $&$    $&$f_{2} $&$+f_{3} $&$ $&$ -e_{4}$ \\
	  $ L_{1001}:= $&$   $&$  $&$f_{3} $&$+f_{4} $&\\
	  $ L_{1010}:= $&$  f_{1}  $&$+f_{2} $&$ +f_{3}$&$ -e_{3}$& \\
	  $ L_{1011}:= $&$ f_{1}  $&$  $&$ $&$ +f_{4}$&\\
	  $ L_{1100}:= $&$ f_{1} $&$ +f_{2} $&$+f_{3} $&$ $&$ -e_{4}$\\
	  $ L_{1101}:= $&$  f_{1} $&$ +f_{2} $&$ $&$ $&\\
	  $ L_{1110}:= $&$  $&$  $&$ $&$l $& \\
	  $ L_{1111}:= $&$ f_{1} $&$  $&$+f_{3}.$&$ $&
	  \end{tabular}
	  $$
	  These divisors $ D_{\sigma}, L_{\chi} $ satisfy the following relations:
	  $$
	  \begin{adjustbox}{max width=\textwidth}
	  \begin{tabular}{l l l l l l l l l l l l r r r r r r}
	  $ 2L_{0001} $& $ \equiv $&$ D_{0011} $& $ +D_{0101} $ & $  $ &$ +D_{0111} $& $ +D_{1001} $& $  $ &$ +D_{1011} $& $ +D_{1101} $& $  $ &$ +D_{1111} $&$ \equiv $&$ 4f_{1} $&$+2f_{2} $&$ $&$-2e_{3} $&\\
	  $ 2L_{0010} $& $ \equiv $&$ D_{0011} $& $  $ & $  +D_{0110} $ &$ +D_{0111} $& $  $& $ +D_{1010} $ &$ +D_{1011} $& $  $& $ +D_{1110} $& $ +D_{1111} $&$ \equiv $&$  $&$ 2f_{2} $&$ $&$+2l $&\\
	  $ 2L_{0100} $& $ \equiv $&& $ D_{0101} $ & $ + D_{0110} $ &$ +D_{0111} $& $  $& $  $ &$  $& $ +D_{1101} $& $ +D_{1110} $& $ +D_{1111} $&$ \equiv $&$ 2f_{1} $&$ +2f_{2} $&$+2f_{3} $&$ $&$ -2e_{4}$\\
	  $ 2L_{1000} $& $ \equiv $&& $  $ & $  $ &$  $& $ D_{1001} $& $ +D_{1010} $ &$ +D_{1011} $& $ +D_{1101} $& $ +D_{1110} $& $ +D_{1111} $&$ \equiv $&$ 2f_{1} $&$ +2f_{2} $&$+2f_{3} $&$ $&$ -2e_{4}$\\
	  $ 2L_{0011} $& $ \equiv $&& $ D_{0101} $ & $ + D_{0110} $ &$  $& $ +D_{1001} $& $ +D_{1010} $ &$  $& $ +D_{1101} $& $ +D_{1110} $& $  $&$ \equiv $&$ 2f_{1}$&$+4f_{2} $&$ $&$-2e_{3} $&$ -2e_{4}$\\
	  $ 2L_{0101} $& $ \equiv $&$ D_{0011} $& $  $ & $  +D_{0110} $ &$  $& $ +D_{1001} $& $  $ &$ +D_{1011} $& $  $& $ +D_{1110} $& $  $&$ \equiv $&$   $&$2f_{2} $&$ +2f_{3}$&$ $&\\
	  $ 2L_{0110} $& $ \equiv $&$ D_{0011} $& $ +D_{0101} $ & $  $ &$  $& $  $& $ +D_{1010} $ &$ +D_{1011} $& $ +D_{1101} $& $  $& $  $&$ \equiv $&$ 4f_{1}  $&$ +2f_{2} $&$ $&$-2e_{3}  $&$ -2e_{4}$\\
	  $ 2L_{0111} $& $ \equiv $&& $  $ & $  $ &$ D_{0111} $& $ +D_{1001} $& $ +D_{1010} $ &$  $& $  $& $  $& $ +D_{1111} $&$ \equiv $&$    $&$2f_{2} $&$+2f_{3} $&$ $&$ -2e_{4}$\\
	  $ 2L_{1001} $& $ \equiv $&$ D_{0011} $& $ +D_{0101} $ & $  $ &$ +D_{0111} $& $  $& $ +D_{1010} $ &$  $& $  $& $ +D_{1110} $& $  $&$ \equiv $&$   $&$  $&$2f_{3} $&$+2f_{4} $&\\
	  $ 2L_{1010} $& $ \equiv $&$ D_{0011} $& $  $ & $  +D_{0110} $ &$ +D_{0111} $& $ +D_{1001} $& $  $ &$  $& $ +D_{1101} $& $  $& $  $&$ \equiv $&$  2f_{1}  $&$+2f_{2} $&$ +2f_{3}$&$ -2e_{3}$&\\
	  $ 2L_{1011} $& $ \equiv $&& $ D_{0101} $ & $ + D_{0110} $ &$  $& $  $& $  $ &$ +D_{1011} $& $  $& $  $& $ +D_{1111} $&$ \equiv $&$ 2f_{1}  $&$  $&$ $&$ +2f_{4}$&\\
	  $ 2L_{1100} $& $ \equiv $&& $ D_{0101} $ & $ + D_{0110} $ &$ +D_{0111} $& $ +D_{1001} $& $ +D_{1010} $ &$ +D_{1011} $& $  $& $  $& $  $&$ \equiv $&$ 2f_{1} $&$ +2f_{2} $&$+2f_{3} $&$ $&$ -2e_{4}$\\
	  $ 2L_{1101} $& $ \equiv $&$ D_{0011} $& $  $ & $  +D_{0110} $ &$  $& $  $& $ +D_{1010} $ &$  $& $ +D_{1101} $& $  $& $ +D_{1111} $&$ \equiv $&$  2f_{1} $&$ +2f_{2} $&$ $&$ $&\\
	  $ 2L_{1110} $& $ \equiv $&$ D_{0011} $& $ +D_{0101} $ & $  $ &$  $& $ +D_{1001} $& $  $ &$  $& $  $& $ +D_{1110} $& $ +D_{1111} $&$ \equiv $&$  $&$  $&$ $&$2l $&\\
	  $ 2L_{1111} $& $ \equiv $&& $  $ & $  $ &$ D_{0111} $& $  $& $  $ &$ +D_{1011} $& $ +D_{1101} $& $ +D_{1110} $& $  $&$ \equiv $&$ 2f_{1} $&$  $&$+2f_{3}.$&$ $&
	  \end{tabular}
	  \end{adjustbox}
	  $$
	  
	  \noindent
	  Thus by Proposition \ref{Construction of cover of degree 8}, the divisors $ D_{\sigma}, L_{\chi} $ define a $ \mathbb{Z}^4_2 $-cover $ \xymatrix{g: X \ar[r] & Y_4}  $. Moreover, this $ \mathbb{Z}^4_2 $-cover fulfils the hypotheses of Theorem \ref{the theorem with d = 20}. In fact, we have 
	  \begin{align*}
	  	D_{0100} + D_{0101} +D_{0110} + D_{0111} &=h_{14} + f_{21}+ f_{31} \equiv 3l - e_1 - e_2 - e_3 - e_4\\ 
	  	D_{1000} + D_{1001} +D_{1010} + D_{1011} &=e_{2} + h_{23} + h_{24} + f_{11} \equiv 3l - e_1 - e_2 - e_3 - e_4\\
	  	D_{1100} + D_{1101} +D_{1110} + D_{1111} &=h_{34} + h_{12} +h_{13} + e_{1} + e_{3} \equiv 3l - e_1 - e_2 - e_3 - e_4,
	  \end{align*}
	  \noindent
	  $ h^{0}\left( K_{Y_4} + L_{\chi} \right) =0$ for all $ \chi \notin \left\lbrace \chi_{1000}, \chi_{0100}, \chi_{1100}\right\rbrace  $, and the divisor $ D_{0001}+ D_{0010}+ D_{0011} - K_{Y_4} \equiv 3l - e_1 - e_2 - e_3$ is nef and big. Thus by Theorem \ref{the theorem with d = 20} and Proposition \ref{invariants of Z_2^n cover}, the surface $ X $ is a minimal surface of general type and possesses the following invariants:
	  \begin{align*}
	  	K_S^2= 24, p_g\left( S\right) = 3, \chi\left( \mathcal{O}_S\right) =4, q\left( S\right)  = 0.  
	  \end{align*}
	  \noindent
	  Moreover, the canonical map $ \varphi_{\left| K_X  \right|}  $ is of degree $ 20 $ and the two $ \left( -2\right) $-curves coming from $ \overline{e}_{4} $ are the fixed part of $ \left| K_X\right|  $. Therefore, we obtain the surface in the second row of Theorem \ref{the main theorem with d = 20}.\\
	  
	  \begin{Remark}
	  The surface $ X $ has three pencils of genus $ 9 $ corresponding the fibres $ f_1, f_2, f_3$ and a pencil of genus $ 13 $ corresponding to the fibre $ f_4$.
	  \end{Remark}
     	  
	  	\begin{Remark}
	  	The surface $ X $ admits natural deformations. Moreover, all the natural deformations of $ X $ are Galois.
	  \end{Remark}
	  \noindent
	  Similarly to Remark \ref{deformation of surface}, we have that $ H^{0}\left( Y_4, D_{0110}\right) \cong H^{0}\left( Y_4, D_{0111}\right) \cong H^{0}\left( Y_4, D_{1011}\right) \cong \mathbb{C}^2$ and $ H^{0}\left( Y_4, D_{\sigma}\right) \cong \mathbb{C} $ for the other non-trivial $ D_{\sigma} $. This implies that the famify of natural deformations of $\xymatrix{g: X \ar[r] & Y_4} $ is parametrized by the base space $ \mathbb{P}^1\times \mathbb{P}^1\times\mathbb{P}^1$. Furthermore, all natural deformations of $ X $ are Galois since $ \bigoplus \limits_{\substack{\sigma \neq 0 \\ \chi \ne \chi_{0000} \\\chi\left( \sigma\right) =1}  }{H^{0}\left( Y_4, D_{\sigma} - L_{\chi}\right) } =0$.

\section*{Acknowledgments}
The author is deeply indebted to Margarida Mendes Lopes for all her help. Thanks are also due to Jungkai Alfred Chen for the suggestion on the existence of deformations. The author would like to express his gratitude to the anonymous referee for his/her thorough reading of the paper and valuable suggestion.

The author was partially supported by Funda\c{c}\~{a}o para a Ci\^{e}ncia e Tecnologia (FCT), Portugal through the program Lisbon Mathematics PhD (LisMath) of the University of Lisbon, scholarship FCT - PD/BD/113632/2015 and project UID/MAT/04459/2019 of CAMGSD. This paper was finished during the author's postdoctoral fellowship at the National Center for Theoretical Sciences (NCTS), Taiwan, under the grant number MOST 109-2119-M-002-014. The author would like to thank NCTS for the financial support and kind hospitality.


\Addresses

\end{document}